\documentclass[11pt]{article}
\usepackage[english]{babel}
\usepackage{amssymb,amsmath,amsfonts}
\usepackage{float}
\usepackage{srcltx}
\usepackage{theorem}
\usepackage{graphicx}
\usepackage{xcolor}
\usepackage{caption}

\date{}
\newcommand{\Csharp}{%
  {\settoheight{\dimen0}{C}C\kern-.05em \resizebox{!}{\dimen0}{\raisebox{\depth}{\#}}}}

\newcommand{\alert}[2][magenta]%
{{\color{#1}\mbox{[\hspace{-0.4ex}[}#2\mbox{]\hspace{-0.4ex}]}}}

\newcommand{\btheorem}{\begin{theorem}}
\newcommand{\etheorem}{\end{theorem}}

\theoremstyle{remark}
\newtheorem{remark}{Remark}

\newtheorem{construction}{Construction}

\newtheorem{theorem}{Theorem}

\newtheorem{corollary}{Corollary}

\newtheorem{proposition}{\bf Proposition}

\newcommand {\bproof }{{\par\medskip\noindent \bf Proof. }}
\newcommand {\eproof }{\hfill $\blacktriangle$ \\ \medskip}

\makeatletter
\renewcommand*{\@fnsymbol}[1]{\ifcase#1\or*\else\@arabic{\numexpr#1-1\relax}\fi}
\makeatother

\def\p0{\parindent 0pt}
\author{Alexey D. Marin\thanks{Alexey D. Marin is with Novosibirsk State University, Novobisibirsk, Russia, a.marin@g.nsu.ru} \and Ivan Yu. Mogilnykh\thanks{Ivan Yu. Mogilnykh is with  Sobolev Institute of Mathematics, Novosibirsk, Russia, ivmog84@gmail.com} }
\title{{\bf \large Binary codes from subset inclusion matrices}\thanks{This study was performed according to the Government research assignment
for the Sobolev Institute of Mathematics, Siberian Branch of the Russian
Academy of Sciences, project FWNF-2022-0017}}

\begin{document}
\maketitle
\begin{abstract}
 In this paper, we study the minimum distances of binary linear codes with parity check matrices formed from subset inclusion matrices 
$W_{t,n,k}$, representing 
$t$-element subsets versus 
$k$-element subsets of an 
$n$-element set. We provide both lower and upper bounds on the minimum distances of these codes and determine the exact values for any 
$t\leq 3$ and sufficiently large $n$. Our study combines design and integer linear programming techniques. The codes we consider are connected to locally recoverable codes, LDPC codes and combinato\-rial designs.

\end{abstract}
\section{Introduction}
Let $W_{t,n,k}$ be the inclusion matrix of $t$-subsets versus $k$-subsets of an $n$-element set: $W_{t,n,k}(U,V)=1$ if $U$ is contained in $V$ and 0 otherwise. Consider the binary linear code $C_{t,n,k}$ with parity check matrix $W_{t,n,k}$: $$C_{t,n,k}=\{x: W_{t,n,k}x^T=0 \}.$$

In \cite{Wilrank} Wilson  obtained the following expression for the rank of $W_{t,n,k}$ (thus solving the dimension problem for the codes $C_{t,n,k}$) over the finite field of prime order $p$:
\begin{equation}\label{WilsonRank} \sum_{i:(^{k-i}_{t-i})\, mod \, p \neq 0} (^n_i)-(^n_{i-1}).\end{equation}

Throughout what follows we denote by $d_{t,n,k}$ the minimum distance of the code $C_{t,n,k}$.
In general, for a fixed $t$ the minimum distances of the codes from Wilson inclusion matrices vary in a wide range:
\begin{equation}\label{general_bounds}
t+2\leq d_{t,n,k} \leq 2^{t+1}.
\end{equation}
Both bounds follow by converting several previous results \cite{WZ}, \cite{Pot}, \cite{BI}[Theorem 2] (lower bound and its sharpness), \cite{Hwang} (upper bound) to terms of the current work, but we provide a self-contained representation in our paper  (see Sections \ref{S_const} and  \ref{S_LB}).  

For $t=k-1$ the matrices $W_{k-1,n,k}$   are free of $2\times 2$ all-ones submatrices and the corresponding Tanner graphs of the codes $C_{k-1,n,k}$ have girth $6$. In \cite{WZ}  these codes are treated as multiuser locally recoverable codes.  Due to a wide range of column and row sum in the Wilson matrices $W_{k-1,n,k}$, the class contains the codes with arbitrary locality and availability, whereas having a better code rate $(^{n-1}_{k})/(^n_k)=1-\frac{k}{n}$ than several  other mass-type constructions, such as a direct product \cite{TB}. The explicit expression $(^{n-1}_{k})$ for the dimension of the code $C_{k-1,n,k}$ follows from the Wilson rank formula (\ref{WilsonRank}). The minimum distances in this particular case were also obtained:

\begin{theorem}\label{T_1}\cite{WZ}, \cite{Pot}
The minimum distance $d_{k-1,n,k}$ of the code $C_{k-1,n,k}$ is $k+1$.
\end{theorem}

This result was independently obtained by Potapov \cite{Pot}, who studied unitrades, a combinatorial equivalent of the nonzero codewords of the code $C_{k-1,n,k}$. The small weight distribution for these codes in \cite{Pot}[Section 3] were shown to have gaps similar to those of Reed-Muller codes with no codewords of weights greater than $k+1$ (minimum distance) and less than $2k$;  the minimum and pre-minimum weight, i.e. $2k$, codewords were also characterized.

By replacing $k$-subset with a $q$-ary $k$-subspace in the definition of the inclusion matrix $W_{t,n,k}$ one obtains its $q$-ary generalization, $W_{t,n,k}^q$. The problem of determining the minimum distance of the $q$-ary codes with these parity check matrices as well as their ranks are  well-known  problems (the latter was solved only for $t=1$ by Hamada \cite{Ham}), which are open in general. Even for the codes with parity check matrices $W_{1,n,k}^q$ the minimum distance is not known in general. 
We refer to the work \cite{L} of Lavrauw et. al and its introductory part for the review of the minimum distance problem for $q$-ary generalization of the codes $C_{t,k,n}$. We note that the authors of \cite{L}  consider the codes from transposed matrix $(W^q_{t,k,n})^T$, $t<k$ which is equivalent to the case $W^q_{t,k,n}$  up to duality for $t=n-t$ and $k=n-k$. For $q=2$ the binary codes with these parity check matrices are known to produce good LDPC codes \cite{KLF}. 

The codes from Wilson-type matrices are conceptually close to other codes arising from natural  incidence structures such as the incidence matrix of Hamming graph  \cite{F}, \cite{F2}. A natural description for the minimum weight codewords is provided in these works and the  permutation decoding is suggested. Despite the fact that (at least partial) permutation decoding for the codes $C_{t,n,k}$ is also possible, we provide a view on these codes (after circulant lifting) as LDPC codes  and focus on  bit flipping and min-sum decoders in Section \ref{S_QC}.

The main problem of the current contribution is determining the value of the minimum distance $d_{t,n,k}$ of the codes $C_{t,n,k}$. We provide a complete answer to this question for all $t\leq 3$ and all $n\geq n_0$, where $n_0$ linearly depends on $k$. 

In the current study of the small weight codewords of $C_{t,n,k}$  we use equivalent terms of binary designs, which are a subcase of $p$-ary designs introduced in \cite{Wil:2009}, rather than codeword terms. For an $n$-element set $N$ (called {\it points}), by a {\it binary} $t$-$(n,k)$-{\it design} (or simply $t$-$(n,k)$-{\it design}) we mean a nonempty set $D$ of $k$-subsets of $N$ (called {\it blocks}) such that any $t$-subset is contained in an even number of blocks of $D$.  We note that any binary $t$-$(n,k)$-design is a  $t$-$(n',k)$-design for all $n'\geq n$.   The binary  $t$-$(n,k)$-designs are nonempty and do not have multifold blocks and therefore are in one-to-one correspondence with the linearly dependent subsets of the columns of the matrix $W_{t,n,k}$ and we have the following.

\begin{proposition}

The nonzero codewords of  $C_{t,n,k}$ are in a one-to-one correspondence with the $t$-$(n,k)$-designs, where the weight of a codeword equals the number of blocks of the corresponding binary design.
\end{proposition}

In Section \ref{S2} we provide necessary notations and definitions for our study of binary designs and  quasi-cyclic LDPC codes from Wilson matrices. Section \ref{S_const} is devoted to  constructions for binary designs with a small number of blocks.  

In Section \ref{S_LB} we exploit classic design approaches for obtaining lower bounds on the number of blocks in binary designs. We note that a binary $t$-design is not necessarily a binary $i$-design for $0\leq i\leq t$, however there are infinite series of linear codes $C_{t,n,k}$ composed of solely of such binary designs,  providing sharpness for the upper bound $2^{t+1}$ on $d_{t,n,k}$ in (\ref{general_bounds}).

In Sections \ref{S_t2}-\ref{S:t=3evenk} we find  the minimum number of blocks in binary $2$-designs and $3$-designs. We split the treatment of  $3$-$(n,k)$-designs into two different subcases by parity of $k$. The case of odd $k$ is handled via lower bounds developed in Section \ref{S_LB},  constructions from Section \ref{S_const} and additional combinatorial considerations. The minimization problem for the number of blocks in $t$-$(n,k)$-designs is solved by a developed integer linear programming search targeted at the reduced incidence matrices of these designs.   

 In Section \ref{S_QC} we discuss quasicyclic LDPC codes from Wilson inclusion matrices and decoding of these codes. It is well-known that the presence of short cycles in the Tanner graph  plays a significant role in degrading LDPC decoding performance. By excluding these cycles via circulant liftings, we obtain  LDPC codes from $W_{t,n,k}$. The obtained codes provide  similar decoding performance  than that of random-type Mackay codes \cite{MacKay} for bit-flipping \cite{GDBF} and layered min-sum \cite{minsum} decoders.

\section{Notations and definitions}\label{S2}

\subsection{Linear codes}
{\it The Hamming distance} between binary tuples is defined as the number of positions in which tuples are different. {\it The weight} of a binary tuple is the number of its ones. 
By a binary {\it code} we mean a subset of binary tuples of length $n$, where the latter is called {\it the length} of the code. The {\it minimum distance} is defined as the minimum Hamming distance between pairs of distinct codewords in the code.  Although the maximum likelihood decoder can correct up to $\left\lfloor\frac{d-1}{2}\right\rfloor$ errors, its practical implementation becomes infeasible even for relatively short code lengths due to severe complexity issues.

A code is {\it linear} if its codewords form a vector space with respect to addition via modulo two.  
The {\it dimension} of a linear code refers to its dimension as a subspace. Notably, the minimum distance of a linear code is equal to the minimum nonzero weight of a codeword in the linear code.  A {\it parity check matrix} of a binary linear code $C$ is a matrix $H$ such that the code is the null space of $H$, that is, $C = \{x : Hx^T = 0\}$. The dimension of a linear code is equal to its length minus the rank (over $GF(2)$) of its parity check matrix. Given an $r \times n$ parity check matrix $H$, one can construct a bipartite graph with vertex set $P\cup V$, $P = {P_1, \ldots, P_r}$ (parity nodes) and $V = {V_1, \ldots, V_n}$ (variable nodes), where an edge exists between $P_i$ and $V_j$ if and only if $H_{i,j} = 1$. This graph is referred to as the {\it Tanner graph}.

LDPC codes represent a broadly defined class of codes characterized by having sparse parity check matrices, known for their excellent decoding performance under iterative algorithms. The decoding  of an  LDPC code is influenced by several factors related to both the code structure and the properties of its associated Tanner graph, see for example \cite{RU}:

$\bullet$ Sparsity of the parity check matrix

$\bullet$ Minimum distance of the code

$\bullet$ Vertex degree distribution in the Tanner graph

$\bullet$ Girth and the presence of short cycles in the Tanner graph

$\bullet$ Various pseudo-codeword structures, including trapping, stopping, and absorbing sets in the Tanner graph.

Wilson-type matrices $W_{t,n,k}$ are sparse and regular, with constant row and column sums. Codes with these matrices as parity check matrices can correct multiple errors, as we will see in Sections \ref{S:t=3oddk} and \ref{S:t=3evenk} where we will find the minimum distances of the codes $C_{t,n,k}$ for $t=1$, $2$, $3$. However, except for the case when $t = k-1$, the Tanner graphs of the codes $C_{t,n,k}$ exhibit numerous cycles of length four and six, which significantly degrade decoding performance. To address this issue, we employ circulant liftings in Section \ref{S_QC}, a technique often used for avoiding short cycles. It is also known to increase the minimum distance to some extent. 

Traditionally, a parity check matrix $H$ of a quasi-cyclic LDPC code is derived from a binary matrix $B$ (referred to as the  {\it base matrix} or {\it mother matrix}), which has a relatively small number of rows and columns. Each zero element of $B$ is replaced by an all-zero matrix of size $qc$, and each unit element $(i, j)$ of $B$ is replaced by a permutation matrix of size $qc$ shifted by $\epsilon_{i,j}$ positions, where $0 \leq \epsilon_{i,j} \leq qc-1$.

The matrix composed of elements $\epsilon_{i,j}$ (with $\epsilon_{i,j} = -1$ if $M_{i,j} = 0$) is called the {\it exponent matrix}. It is well-known that the properties of the base code $\{x : Bx^T = 0\}$, including its minimum distance,  influence those of the quasi-cyclic code with the parity check matrix $H$.

 It is obvious that any 
 cycle in Tanner graph of base code gives rise to either zero or $qc$ cycles in the Tanner graph of the quasi-cyclic code. Moreover, the absence of cycles of length four  in the Tanner graph of a quasi-cyclic code with exponent matrix $\epsilon$ can be expressed using simple condition in terms of the elements of exponent matrix $\epsilon$.
For any row $a$, $b$ and columns $a'$, $b'$ such that $B_{a,a'}=B_{a,b'}=B_{b,a'}=B_{b,b'}=1$ (a $4$-cycle in the base code) 
we have:
\begin{equation}\label{C4}   \epsilon_{a,a'}+  \epsilon_{b,b'}\neq  \epsilon_{a,b'}+  \epsilon_{b,a'} \mbox{ mod }qc.\end{equation}
Similarly, the property of Tanner graph being  cycle six free is formulated as follows.
For any rows $a$, $b$, $c$ and columns $a'$, $b'$, $c'$ such that the $B_{a,a'}=B_{b,b'}= B_{c,c'}=  B_{a,b'}=B_{b,c'}=B_{c,a'}=1$ (a $6$-cycle in the base code)  we have
\begin{equation}\label{C6}   \epsilon_{a,a'}+  \epsilon_{b,b'}+ \epsilon_{c,c'}\neq  \epsilon_{a,b'}+\epsilon_{b,c'}+\epsilon_{c,a'} \mbox{ mod }qc.\end{equation}

We finish with an example demonstrating the notions and objects above.

{\it Example 1.}
Let the exponent matrix of the LDPC code be 
$$\left(%
\begin{array}{cccc} 2& 1 & 0 &1 \\
0& 1 & -1 &0 \\
-1& 0 & 2 &1\\\end{array}%
\right),$$ the size of circulant be $qc=3$, the base matrix is as follows
$$B=\left(%
\begin{array}{cccc} 1& 1 & 1 &1 \\
1& 1 & 0 &1 \\
0& 1 & 1 &1\\\end{array}%
\right).$$
The corresponding parity check matrix of the binary quasicyclic LDPC code is 

$$\left(%
\begin{array}{ccc|ccc|ccc|ccc}
0 & 0 & 1 & 0 & 1 & 0 & 1 & 0 & 0 & 0 & 1 & 0 \\
1 & 0 & 0 & 0 & 0 & 1 & 0 & 1 & 0 & 0 & 0 & 1 \\
0 & 1 & 0 & 1 & 0 & 0 & 0 & 0 & 1 & 1 & 0 & 0 \\\hline
1 & 0 & 0 & 0 & 1 & 0 & 0 & 0 & 0 & 1 & 0 & 0 \\
0 & 1 & 0 & 0 & 0 & 1 & 0 & 0 & 0 & 0 & 1 & 0 \\
0 & 0 & 1 & 1 & 0 & 0 & 0 & 0 & 0 & 0 & 0 & 1 \\\hline
0 & 0 & 0 & 1 & 0 & 0 & 0 & 0 & 1 & 0 &1 & 0 \\
0 & 0 & 0 & 0 & 1 & 0 & 1 & 0 & 0 & 0 & 0 & 1 \\
0 & 0 & 0 & 0 & 0 & 1 & 0 & 1 & 0 & 1 & 0 & 0 \\
\end{array}%
\right).$$ 
We see that the Tanner graph of the code contains $3$ cycles of length $4$ and $3$ cycles of length $6$. This can be viewed as respective paths in the parity check matrix as well as from equalities (\ref{C4}) and (\ref{C6}): 
$$\epsilon_{1,2}+\epsilon_{3,3}=\epsilon_{1,3}+\epsilon_{3,2} \, \mbox{ mod }3,$$
$$\epsilon_{1,1}+\epsilon_{3,2}+\epsilon_{2,4}=\epsilon_{1,2}+\epsilon_{3,4}+\epsilon_{2,1} \, \mbox{ mod }3.$$
Using MAGMA \cite{MAGMA}  we found that the minimum distance of the lifted LDPC code is $4$ which is improved compared to that of the base code $\{x:Bx^T=0\}$, which is $2$.
\subsection{Designs}

For a subset $D$ of $k$-subsets and a set $S$ we denote by $\lambda(S,D)$ the number of  subsets of $D$ containing $S$. The definition of a binary $t$-$(n,k)$-design $D$ ($p$-ary $t$-$(n,k,\lambda)$ design \cite{Wil:2009} and classic $t$-$(n,k,\lambda)$-design respectively) is equivalent to  $\lambda(S,D)=0\mbox{ mod  2}$ ($\lambda(S,D)=\lambda\mbox{ mod }  p$  and $\lambda(S,D)=\lambda$ respectively) for all subsets $S$ of $t$ points. We note that our definition of binary design is a subcase of $p$-ary designs for $p=2$ introduced by Wilson in \cite{Wil:2009}. The $p$-ary designs were also considered in  \cite{Liam}.

For a set $D$ of $k$-subsets and a point $i$, we introduce the following notation:
$$D^{i}=\{B:B \in D, i\in B\},$$
$$(D^{i})'=\{B\setminus i:B \in D, i\in B\}.$$

For a classical design $D$, the set $(D^{i})'$ defined above is called the {\it derivative design} of $D$ with respect to the point $i$. We note the following obvious equality $\lambda(i,D)=|D^i|=|(D^{i})'|$.

A point-block {\it  incidence matrix} of a set  $D$ of subsets of the point set $N$ is defined as a matrix $A$, whose rows are indexed by the elements of the set $N$, and columns by the subsets from $D$. The entry $A_{i,B} = 1$ if and only if the element $i$ belongs to the subset $B$ from $D$, and $A_{i,B} = 0$ otherwise. The \textit{reduced incidence matrix} of a binary $t$-$(n,k)$-design $D$ is the matrix $A$ obtained from the incidence matrix by removing duplicate rows.

At times, the property of being a binary $t$-$(n,k)$-design can be  more clearly represented by the following relation of the rows of its incidence matrix: any $t$ rows (with distinct numbers) of the incidence matrix of such a design have an even number of common ones. We proceed the following simple example.

{\bf Example 2.} Let $D$ be $\{\{1,3,4\},\{1,2,4\},\{2,3,4\},\{1,3,5\},\{1,2,5\},\{2,3,5\}\}$ with the incidence matrix  $$\left(%
\begin{array}{cccccc} 1& 1 & 0 &1 &1&0\\0& 1 & 1 &0 &1&1\\
1& 0 & 1 &1 &0&1\\1& 1 & 1 &0 &0&0\\ 0&0&0&1&1&1\end{array}%
\right).$$ As can be easily seen, any pair of rows with different indices have an even number of common ones. We conclude that 
$D$ is a $2$-$(6,3)$-design.

\section{Constructions of binary designs}\label{S_const}

\begin{construction}\label{Constr_1} (Subsets).\end{construction} Let $D$ be the set of all $k$-subsets of a $(k+1)$-set. If $k-t$ is odd, then $D$ is a binary $t$-$(k+1,k)$-design. Indeed, any $t$-subset is contained in $k+1-t$ (an even number) subsets from $D$.

Therefore we obtain the following bound for all $n, \,n\geq k+1$ and odd $k-t$: 
\begin{equation}\label{Ub_constr_1}
 d_{t,n,k}\leq k+1
\end{equation}

\begin{construction}\label{Constr_2}
(Doubling). \end{construction} Let $D$ be a $t$-$(n,k)$-design. Consider two additional points $n+1$ and $n+2$ and the  block set $\{B \cup \{n+1\} : B \in D\} \cup \{B \cup \{n+2\} : B \in D\}$, which we denote by $D^{+}$. We show that $D^{+}$ is a $(t+1)$-$(n+2,k+1)$-design.

Let $S$ be a set of $t+1$ points. If $S$ contains both $n+1$ and $n+2$, then it is not contained in any block from $D^{+}$ by the construction of this set.
If $S$ does not contain $n+1$ and $n+2$, then $S$ is contained in $\lambda(S,D)$ blocks of $D$.
Since the blocks of $D^{+}$ are obtained by adding points $n+1$ and $n+2$ to the blocks of $D$, the set
$S$ is contained in exactly $2\lambda(S,D)$ blocks from $D^{+}$. If $S$ contains $n+1$ but not $n+2$, then, by definition of $D^{+}$, we have $\lambda(S,D^+) = \lambda(S,D)$, which is even since $D$ is a $t$-$(n,k)$-design. Hence, we have the following bound:
\begin{equation}\label{ubpund1}
d_{t+1,n+2,k+1} \leq 2d_{t,n,k}.
\end{equation}

We note that the design from Example 2 is obtained by applying the doubling construction to a $1$-$(3,2)$ design with blocks $\{1,2\}$, $\{2,3\}$, $\{1,3\}$, which in turn is obtained as Construction \ref{Constr_1}.

\begin{construction}\label{Constr_3}
(Generalized Pasch-configuration). \end{construction}

The point set of this $t$-$(k+t+1,k)$ design is $\{1,\ldots,k+t+1\}$. For each point $i\in \{1,\ldots,t+1\}$, denote by $f(i)$ the point $i+t+1$, and for a subset $R\subseteq \{1,\ldots,t+1\}$, let  $f(R)$ denote $\{f(r) : r \in R\}$. Let  $D$ be 
$$\{R\cup f(\{1,\ldots,t+1\}\setminus R)\cup \{2t+3,\ldots, k+t+1\}: R\subseteq \{1,\ldots,t+1\}\}.$$
We see that $|D|=2^{t+1}$.
In \cite{Hwang} it is shown that the set of blocks $D$ above can be partitioned into two subsets $T$ and $T'$ such that for any $t$-set is contained in the same number of blocks from $T$ and $T'$. It follows that $D$ is a binary $t$-$(k+t+1,k)$-design.

\begin{proposition}\label{P_H} \cite{Hwang}
For all $n \geq k+t+1$ the following holds:
$$d_{t,n,k}\leq 2^{t+1}.$$
\end{proposition}

\begin{construction}\label{Constr_4}  (Binary designs from the  $2$-$(7,4,2)$ Hadamard  design).\end{construction}  Consider the blocks of the classical $2$-$(7,4,2)$ Hadamard design $D_H$:
$$D_H=\{\{1,2,3,4\},\{1,2,5,6\},\{3,4,5,6\},\{1,4,6,7\},\{1,3,5,7\},\{2,4,5,7\},\{2,3,6,7\}\}.$$ It is obvious that the classical $2$-$(7,4,2)$ Hadamard design $D_H$ is a binary $2$-$(7,4)$-design. Moreover, since for any $i \in \{1,\ldots,7\}$ we have $\lambda(i, D_H) = 4$, the Hadamard design is also a $1$-$(7,4)$-design. The Hadamard design $D_H$ is not a binary $0$-design as it contains an odd number of blocks.

Let $l$ be an arbitrary integer, $l \geq 1$.
We increase the block size by a factor of $l$.
To achieve this, we replace each block $B$ of the design $D_H$ with the block
$\{i+7j:i\in B,j\in \{0,\ldots,l-1\}\}$. We denote the resulting set of seven blocks of size $4l$ by $D$. If two elements $a, b \in \{1, \ldots, 7l\}$ are such that $a-b$ is divisible by $7$, then $\lambda(\{a,b\}, D) = \lambda((a \mod 7) + 1, D_H) = 4$. Otherwise, $\lambda(\{a,b\}, D) = \lambda(\{(a \mod 7) + 1, (b \mod 7) + 1\}, D_H) = 2$. We conclude that $D$ is a binary $2$-$(7l, 4l)$-design with $7$ blocks.
 In particular, for $k$ divisible by $4$, we have $d_{2,n,k} \leq 7$.
 
\begin{remark}\label{remark_1}
The above construction implies that  lower bounds on the number of blocks, similar to Fisher's inequality for binary $2$-designs do not hold, unlike for classical $2$-designs, and the number of points can be arbitrarily large for a fixed number of blocks.
A similar property was previously noted for a more general class of $p$-ary designs in \cite{Wil:2009}.
\end{remark}

\section{Lower bounds on the number of blocks in binary designs}\label{S_LB}
In this section, we provide estimates for the number of blocks in binary designs based on standard approaches for classical designs.

\begin{theorem}\label{T:main}

1. Let $D$ be a $t$-$(n,k)$-design. If $j \in\{0, \ldots, t-1\}$ and $(^{k-j}_{t-j})$ is odd, then $D$ is a $j$-$(n,k)$-design.

2. Let $D$ be a $t$-$(n,k)$-design. For each point $i \in \{1, \ldots, n\}$, the set $(D^{i})'$ is either a $(t-1)$-$(n-1,k-1)$-design or an empty set. In particular, we have
\begin{equation}\label{e:2}
d_{t,n,k} \geq d_{t-1,n-1,k-1} + 1.
\end{equation}

3. Let $D$ be both a $j$-$(n,k)$-design and a $(j-1)$-$(n,k)$-design.  Then for any $i \in \{1, \ldots, n\}$, $D^{i}$ is a $(j-1)$-$(n,k)$-design or empty and $D \setminus D^{i}$ is either a $(j-1)$-$(n-1,k)$-design or an empty set. In particular,
\begin{equation}\label{e:1}
d_{j,n,k} \geq \max\{d_{j-1,n,k}, d_{j-1,n-1,k-1}\} + d_{j-1,n-1,k}.
\end{equation}
If the number of blocks in the design $D$ equals $\max\{d_{j-1,n,k}, d_{j-1,n-1,k-1}\} + d_{j-1,n-1,k}$, then each point in the design $D$ is incident to no blocks of $D$, all blocks of the design $D$ or exactly \,\,\,$\max\{d_{j-1,n,k}, d_{j-1,n-1,k-1}\}$ blocks of the design $D$.

4. If $k-j$ is even then any $j$-$(n,k)$-design $D$ is also a $(j-1)$-$(n,k)$-design. In particular, Theorem \ref{T:main}.3 holds for $D$.
\end{theorem}

\bproof
1. For each $j$-subset  $S$ of points, count 
$|\{(T,B): S\subseteq T\subseteq B, |T|=t, B\in D\}|$ in two ways and obtain the equality 
$$(^{k-j}_{t-j})\lambda(S,D)=\sum\limits_{S\subseteq T, |T|=t}\lambda(T,D).$$ 
Considering the equality modulo $2$ and noting that $\lambda(T,D)$ is even and $(^{k-j}_{t-j})$ is odd, it follows that $\lambda(S,D)$ is even, i.e., $D$ is a $j$-design.

2. For any $(t-1)$-subset $S$ of $\{1, \ldots, n\} \setminus i$,  we have $\lambda(S, (D^i)') = \lambda(S \cup {i}, D)$, which implies that $(D^i)'$ is either empty or a $(t-1)$-$(n-1,k-1)$-design. Clearly, there exists a point $i$ not contained in at least one block of $D$, so $|D| \geq |(D^i)'| + 1$ and the inequality (\ref{e:2}) holds.

3. Let $S$ be an arbitrary $(j-1)$-subset of the  point set. If $S$ does not contain $i$, then it is contained exactly in $\lambda(S \cup {i}, D)$ blocks of $D^i$, which is even because $D$ is a $j$-$(n,k)$-design. Otherwise, $S$ is contained in $\lambda(S, D)$ blocks of $D^i$, which is even because $D$ is a $(j-1)$-$(n,k)$-design. Thus, $D^i$ is a $(j-1)$-design or an empty set. In case both $D^i$ and $D$ are $(j-1)$-$(n,k)$-designs, $D \setminus D^i$ is either a $(j-1)$-$(n-1,k)$-design or an empty set.

Note that due to the proven second statement of the theorem, $(D^i)'$ is a $(j-1)$-$(n-1,k-1)$-design, from which we have  $$|D^i|=|(D^i)'|\geq max\{d_{j-1,n,k},d_{j-1,n-1,k-1}\}$$

Clearly, there is a point $i$ not contained in at least one block of $D$, giving that $D \setminus D^i$ is a non-empty set, so we have:
$$|D|=|D^i|+|D\setminus D^i|\geq max\{d_{j-1,n,k},d_{j-1,n-1,k-1}\}+d_{j-1,n-1,k}. $$
The lower bound above is attained only when the design $D$ is such that for any point $i \in \{1, \ldots, n\}$, we have $\lambda({i}, D) = |D|$ or $\max\{d_{j-1,n,k}, d_{j-1,n-1,k-1}\}$.

4. Since $(^{k-j+1}_{j-j+1}) = k-j+1$ is odd, a $j$-$(n,k)$ design $D$ is a $(j-1)$-$(n,k)$-design by the first statement of the theorem. Thus, by the third statement of the theorem, the required result holds.
\eproof

\begin{remark}\label{remark_2}
Let $i$ be any point of a $t$-$(n,k)$-design $D$ that is incident to a non-zero number of blocks. By Theorem \ref{T:main}.2, the derivative $(D^i)'$ is a $(t-1)$-$(n-1,k-1)$-design, so $\lambda(i,D) \geq d_{t-1,n-1,k-1}$ for blocks of $D$. In other words, any non-zero row of the incidence matrix of $D$ has at least $d_{t-1,n-1,k-1}$ ones.
\end{remark}
A binary $t$-$(n,k)$ design is called {\it universal} \cite{Liam} if it is an $i$-design for all $0 \leq i \leq t$. Due to the following assertion, we see that there exist linear spaces formed by such designs, and the minimum weight in such a space is the largest among all other spaces formed by binary $t$-designs as it attains (\ref{general_bounds}).

\begin{theorem}\label{T2}
Let $k$ be $t \mbox{ modulo } 2^{\lceil log_2(t)\rceil}$. Then for any 
$n>k$, a $t-(n,k)$-design $D$
is a $s$-$(n,k)$-design for all $0\leq s\leq t$. Moreover, for $n\geq k+t+1$ we have at least $2^{t+1}$ blocks in this design and $d_{t,n,k}=2^{t+1}$.

\end{theorem}
\bproof
We show that $(^{k-s}_{t-s})$  is odd for all $s: 0\leq s\leq t$. 

For $l$  we denote by $l_j$ its $j$-th bit in its binary representation: 
$$l=\sum_{j=0,\ldots, \lceil log_2(l)\rceil-1}2^jl_j.$$

By Lucas theorem, we have  
$$(^{k-s}_{t-s})\mbox{ mod 
 }2=\prod\limits_{j=0..\lfloor log_2(t-s)\rfloor} (^{(k-s)_j \mbox{ mod } 2}_{(t-s)_j \mbox{ mod } 2}).$$

Given that $k = t \mod 2^{\lceil \log_2(t) \rceil}$, the least significant bits of the binary representations of $k-s$ and $t-s$ for $s=0,\ldots,t$ are equal: $(k-s)_j = (t-s)_j = 1$, $j=0,\ldots,\lfloor \log_2(t-s) \rfloor$. From this, we obtain that $\binom{k-s}{t-s}$ is odd. By Theorem \ref{T:main}.1, $D$ is an $s$-design for all $0 \leq s \leq t$.

To prove the lower bound on the number of blocks, we use induction and assume that any design that is an $s$-design for all $0 \leq s \leq t$ has at least $2^{t+1}$ blocks. A design $D$ that is an $s$-design for all $0 \leq s \leq t+1$, for a point $i \in \{1,\ldots,n\}$, is partitioned into $D^i$ and $D \setminus D^i$, which are $s$-designs for all $s \in \{0,\ldots,t\}$ according to Theorem \ref{T:main}.3 and thus satisfy the induction hypothesis, having at least $2^{t+1}$ blocks in each.

On the other hand,  by Construction \ref{Constr_3}, for $n \geq k + t + 2$ the minimum number of blocks in a $(t+1)$-$(n,k)$-design is at most $2^{t+2}$ and   the required result is obtained.
\eproof

Note that an arbitrary set of blocks of fixed size is a $0$-design if and only if the number of blocks is even. Thus, for any $k$, $n$ with $n > k \geq 1$, we have $d_{0,n,k} = 2$.

\begin{proposition}\label{p2}

1. If $k$ is odd and $n \geq k + 1$, then $d_{1,n,k} = 4$.

2. If $k$ is even and $n \geq k + 2$, then $d_{1,n,k} = 3$.
\end{proposition}

\bproof By Theorem \ref{T:main}.2, for any $k$ and $n > k$, we have $d_{1,n,k} \geq d_{0,n-1,k-1} + 1$. Since $d_{0,n-1,k-1} = 2$, it follows that $d_{1,n,k} \geq 3$.

For odd $k$, by Theorem \ref{T:main}.4, we have $d_{1,n,k} \geq d_{0,n,k} + d_{0,n-1,k}$. Given that $d_{0,n-1,k} = d_{0,n,k} = 2$ for $k \geq 2$, we obtain the inequality $d_{1,n,k} \geq 4$. For any odd $k$ and $n \geq k + 1$, a $1$-$(n,k)$-design with four blocks is obtained by Construction \ref{Constr_1}. By Construction \ref{Constr_3}, for even $k$, we have $d_{1,n,k} = 3$ for all $n \geq k + 2$.
\eproof
\begin{corollary}\label{C1}
For $n\geq k+t+1$ we have 
 $t+2 \leq d_{t,n,k}\leq 2^{t+1}$.
\end{corollary}
\bproof
The lower bound follows from  Proposition \ref{p2} and application of Theorem \ref{T:main}.2, whereas the upper bound is due to Proposition \ref{P_H}.

\eproof
\section{Minimum number of blocks in a binary $2$-design}\label{S_t2}
\begin{theorem}\label{T:t2} 

1. \cite{WZ}, \cite{Pot} For all $n\geq 4$ we have $d_{2,n,3}=4$.

2. If $k$ is odd, $k\geq 5$ and $n\geq k+3$ then $d_{2,n,k}=6$. Moreover, in every $2$-$(n,k)$-design with $6$ blocks there is a point contained in at least $4$ blocks.

3. If $k=0 \mbox{ mod } 4$, $n\geq k+2$ then we have $d_{2,n,k}\geq 7$. If $n\geq \frac{7k}{4}+2$ then we have $d_{2,n,k}=7$.

4. If $k=2 \mbox{ mod }4$, $n\geq k+1$ then  we have $d_{2,n,k}\geq 8$. If $n\geq k+3$, then $d_{2,n,k}= 8$.

\end{theorem}
\bproof
1. This is a particular case of Theorem \ref{T_1}.

2. Suppose $k$ is odd and not less than $5$. By Theorem \ref{T:main}.2 and Proposition \ref{p2}.2, we have $d_{2,n,k}\geq d_{1,n-1,k-1}+1=4$ for $n\geq k+2$. On the other hand, by the inequality (\ref{ubpund1}) the following holds: 
$$d_{2,n,k}\leq 2d_{1,n-2,k-1}.$$
Hence, given that $n\geq k+3$ by Proposition \ref{p2}.2 we have $d_{1,n-2,k-1}=3$, thus,
$$4\leq d_{2,n,k}\leq 2d_{1,n-2,k-1}=6.$$
Assume there exists a $2$-$(n,k)$-design $D$, $4\leq |D|\leq 6$. We will show contradiction by analyzing the incidence matrix of the design $D$. Since $k\geq 5$, each column of the incidence matrix contains at least $5$ ones. The subsequent considerations are based on the following facts.

$\bullet$ Note that since $d_{1,n-1,k-1}=3$, each nonzero row of the incidence matrix contains at least three ones by Remark \ref{remark_2}.

$\bullet$ Also, each vector of odd weight appears at most once as a row of the incidence matrix of the design, otherwise the incidence matrix of the binary $2$-$(n,k)$-design  contains pairwise non-orthogonal rows.

Let us consider the following cases.

$2a.$ 
Suppose $|D|=4$, i.e. the incidence matrix is a $n\times 4$ matrix. There exists a row of the incidence matrix with weight three, hence there is no row of the incidence matrix with weight four (otherwise, the rows are not pairwise orthogonal). The remaining rows are just  some pairwise distinct rows of weight three. Since there are not more than four such rows, the sum of ones in a column of the incidence matrix is not more than three, which contradicts the condition $k\geq 5$.

$2b'.$ Suppose $|D|=5$ or $6$ and at least $|D|-1$ rows of the incidence matrix have weight $3$. Without loss of generality, since the column sum is $k$, $k\geq 5$ up to permutation of rows the first two rows of the incidence matrix have weight $3$ and at least one common one. Moreover, the first two rows then have exactly $2$ common ones due to pairwise orthogonality of the rows of the incidence matrix. Depending on whether $|D|=5$ or $6$, up to permutation of columns, the incidence matrix is as follows

\begin{equation}\label{M5}\left(%
\begin{array}{ccccc} 1& 1 & 1 &0 &0\\0& 1 & 1 &1 &0\\
.&.&.&.&.\end{array}%
\right)\mbox{ or  }\end{equation}\begin{equation}\label{M6}\left(%
\begin{array}{cccccc} 1& 1 & 1 &0 &0&0\\0& 1 & 1 &1 &0&0\\
.&.&.&.&.&.\end{array}%
\right).
\end{equation}
We proceed with the proof for the matrix (\ref{M6}). The basis of the space of vectors which are orthogonal modulo $2$ to the first two rows of the first matrix consists of the following vectors:
$(110100)$,
$(011000)$,
$(000010)$, $(000001)$. In the space spanned by these vectors, there are only the following three vectors of weight $3$:
$(110100)$,
$(011010)$,  $(011001)$. We see that the first and second, as well as the first and third of these vectors are not orthogonal modulo $2$. Moreover, each vector of weight $3$ can appear not more than once as a row in the incidence matrix, so there are at most $4$ rows of weight $3$ in the incidence matrix. This contradicts that for the current case we have at least $|D|-1=5$ rows of weight three.  
The proof for the matrix (\ref{M5}) (i.e., for $|D|=5$) is obtained from the above by removing all considered vectors in the last position and excluding the vector $(011001)$ from consideration.

$2b''.$ Suppose that $|D| = 5$ or $6$ and all rows of the incidence matrix have weight $3$. Then since $k \geq 5$, there are at least $5$ such rows. From the case $2b'$, we obtain that such an incidence matrix does not exist. 

$2b'''$. 
Suppose that $|D| = 5$ and there are rows of the incidence matrix of weight $4$ or $5$. Obviously, there are no two  distinct  rows orthogonal modulo $2$  of length $5$ with weight $4$. Therefore, in the incidence matrix, there will be a column where ones appear only in rows of weight $3$ and $5$. Since these rows of odd weight can appear at most one time each  and because there are at most $3$  rows of weight $3$  by case $2b'$, the column sum of the incidence matrix is at most $4$, which contradicts $k \geq 5$.

We have shown that $|D| = 4$ or $5$ is impossible. From the fact that case $2b''$ is impossible, it follows that any incidence matrix of a $2$-$(n,k)$ design $D$ with odd $k$ and $6$ blocks contains a row with at least $4$ ones. Equivalently, there is a point incident to at least $4$ blocks of the design $D$.

3. By Theorem \ref{T:main}.4 for $j=2$ and $k$ divisible by $4$ and Proposition \ref{p2},  we have: $$d_{2,n,k}\geq max\{d_{1,n,k},d_{1,n-1,k-1}\}+d_{1,n-1,k}=4+3=7.$$
From Construction \ref{Constr_4}, we obtain that $d_{2,n,k}=7$ for $n\geq \frac{7k}{4}+2$.

4. Follows from Theorem \ref{T2}.
\eproof



\section{Minimum number of blocks in a $3$-$(n,k)$-design for odd $k$}\label{S:t=3oddk}

\begin{proposition}\label{P:1}
For any odd $k$, $k \geq 2$ and $n \geq k+4$, every $3$-$(n,k)$-design contains at least $14$ blocks.  

\end{proposition}
\bproof
Since $k-3$ is even, by Theorem \ref{T:main}.4, a $3$-$(n,k)$-design $D$ is also a $2$-$(n,k)$-design, and we have:
$$d_{3,n,k}\geq max\{d_{2,n,k},d_{2,n-1,k-1}\}+d_{2,n-1,k}.$$
 Since $n \geq k+4$, the conclusion of Theorem \ref{T:t2}.3 and  Theorem \ref{T:t2}.4 holds for $2$-$(n-1,k-1)$-designs, so $d_{2,n-1,k-1}\geq 7$. Consider Theorem \ref{T:t2}.2 for $2$-$(n-1,k)$-design. Since $n \geq k+4$, the condition of Theorem \ref{T:t2}.2 is satisfied and we obtain $d_{2,n-1,k}=6$.
 Hence, we have the inequality:
$$d_{3,n,k}\geq max\{d_{2,n,k},d_{2,n-1,k-1}\}+d_{2,n-1,k}\geq 7+6=13.$$

Assume that $D$ is a $3$-$(n,k)$-design with $13$ blocks. According to Theorem \ref{T:main}.3, each point is incident  to no blocks, all blocks or  exactly to $7$ blocks. The case where there is a point incident to exactly $13$ blocks is not possible. Assume the opposite and such a point $p_1$ exists. Clearly, there must be a point $p_2$ that is not incident to all blocks. By  Theorem \ref{T:main}.3, this point $p_2$ would be incident to exactly $7$ blocks. Thus, the set $\{p_1, p_2\}$ is contained in an odd number of blocks (specifically $7$), which contradicts the fact that $D$ is a $2$-design.

Thus every point in the design $D$ is incident to exactly $7$ blocks or no blocks. Consider an arbitrary point $i$ incident to $7$ blocks. According to Theorem \ref{T:main}.3, $D \setminus D^i$ is a $2$-$(n-1,k)$-design. Given that $\lambda({i}, D) = 7$, we have $|D \setminus D^i| = 13 - 7 = 6$. By Theorem \ref{T:t2}.2 appied to $2$-$(n-1,k)$-design $D \setminus D^i$ with $6$ blocks there must be a point $j \in \{1, \ldots, n\} \setminus \{i\}$ that is incident to at least $4$ blocks of $D \setminus D^i$.

Moreover, since $|D| = 13$ and both $i$ and $j$ are incident to $7$ blocks, there must be at least one block of $D$ containing both points $i$ and $j$. We will show that there are at least $4$ such blocks.

Since $D$ is a $3$-$(n,k)$-design, by Theorem \ref{T:main}.1,  $(D^i)'$ is a $2$-$(n-1,k-1)$-design. Consider the design $\widetilde{D}$ obtained from $(D^i)'$ by taking the derivative with respect to the point $j$. Since both $i$ and $j$ are incident to at least one block of the design $D$, $\widetilde{D}$ is nonempty, so by Theorem \ref{T:main}.1,  $\widetilde{D}$ is a $1$-$(n-2,k-2)$-design and  $|\widetilde{D}|\geq d_{1,n-2,k-2}$. By Proposition \ref{p2} because $k$ is odd we have $|\widetilde{D}| \geq d_{1,n-2,k-2} = 4$.
In other words, by Remark \ref{remark_2}, the points $i$, $j$ are incident to $|\widetilde{D}|\geq 4$ blocks of $D$. 

Thus, $j$ is incident to at least $4$ blocks of $D \setminus D^i$ and  at least $\lambda(j,D^i)=|\widetilde{D}| \geq 4$ blocks of $D^i$, at least $8$ blocks of $D$ total. This contradicts the fact that $j$ is incident to exactly $7$ blocks of $D$, a contradiction.

\eproof
\begin{theorem}
\label{T:t=3}
For any odd $k$  we have
\begin{equation}
d_{3,n,k} = 
 \begin{cases}

   14&\text{  if }k=1 \text{ mod } 4 \text{ and } n\geq max\{k+4,\frac{7k+1}{4}\},
   \\
   16&\text{  if }k=3 \text{ mod } 4 \text{ and } n\geq k+4.
 \end{cases}
\end{equation}

\end{theorem}
\bproof
Consider the cases for $k$.

 Let  $k$ be such that $k=1\mbox{ mod } 4$.  The minimum number of blocks $d_{3,n,k}$ is at least $14$ for $n\geq k+4$, see Proposition \ref{P:1}. 
 On the other hand, since $k=1\mbox{ mod } 4$ one might  consider a $2$-$(\frac{7(k-1)}{4},k-1)$-design from Construction \ref{Constr_4} and further apply Construction \ref{Constr_2} to obtain a $3$-$(\frac{7k+1}{4},k)$-design with $14$ blocks.

The case $k$ such that  $k=3\mbox{ mod } 4$ follows from Theorem \ref{T2}.
\eproof

\section{Minimum number of blocks in a $3$-$(n,k)$-design for even $k$}\label{S:t=3evenk}
\subsection{The properties of reduced incidence matrices}
There are constructions of  binary designs (for example, Construction \ref{Constr_4}) whose incidence matrices can contain a large number of identical rows. For this reason in our search we focus on reduced matrices rather than incidence matrices.

We note that that there are cases where the same matrix can be obtained by reducing several (possible infinitely many)  different incidence matrices. We will illustrate this effect in below.

{\bf Example 2.}
Consider the following matrices $A_l$:

\begin{equation}\label{e1}\left(%
\begin{array}{cccccc} 1& 1 & 1 &0 &0&0\\0& 0 & 0 &1 &1&1\\\hline
1& 1 & 0 &1 &1&0\\
0& 1 & 1 &0 &1&1\\
1& 0 & 1 &1 &0&1\\
\hline
&&.&.\\
&&.&.\\
\hline
1& 1 & 0 &1 &1&0\\
0& 1 & 1 &0 &1&1\\
1& 0 & 1 &1 &0&1\\
\end{array}%
\right),\end{equation}
where the submatrix $\left(%
\begin{array}{cccccc}
1& 1 & 0 &1 &1&0\\
0& 1 & 1 &0 &1&1\\
1& 0 & 1 &1 &0&1\\
\end{array}%
\right)$ is repeated $l$ times.

It is easy to see that any pair of rows with different indices from the matrix $A_l$ are pairwise orthogonal, and the number of ones in a column is constant and equal to $2l+1$. Therefore, the matrix $A_l$ is the incidence matrix of a binary $2$-$(3l+1,2l+1)$-design, which we denote by $D_l$. The reduced incidence matrix of the designs $D_l$ for all $l\geq 1$ is the following matrix: 
$$\left(%
\begin{array}{cccccc}
1& 1 & 1 &0 &0&0\\0& 0 & 0 &1 &1&1\\
1& 1 & 0 &1 &1&0\\
0& 1 & 1 &0 &1&1\\
1& 0 & 1 &1 &0&1\\
\end{array}%
\right).$$

We note the following properties of the incidence matrix of $3$-$(n,k)$-designs and its reduction.

(P1) The columns of the incidence matrix and the reduced incidence matrix are pairwise distinct. The blocks of a binary design are pairwise distinct, so the columns of the incidence matrix are also pairwise distinct. This property is obviously preserved after removing duplicate rows, so it also holds for the reduced incidence matrix.

We will say that a matrix satisfies $3$-\textit{orthogonality} if any three rows with different indices have an even number of columns entirely consisting of ones.

(P2) Both the incidence matrix and the reduced incidence matrix fulfill $3$-orthogonality.

(P3) Each column of the incidence matrix contains $k$ ones. Each column of the reduced incidence matrix contains not more than $k$ ones.

(P4) Each nonzero row of the incidence matrix and the reduced incidence matrix of a $3$-$(n,k)$-design has at least $d_{2,n-1,k-1}$ ones. This follows from  Remark \ref{remark_2}.

(P4$'$) If $k$ is even and not less than $6$, and $n\geq k+2$, then each row of the reduced incidence matrix contains at least $6$ ones.
By Theorem \ref{T:t2}.2, we have $d_{2,n-1,k-1}=6$ for even $k$. This property follows from (P4).

\subsection{Integer linear programming for $3$-$(n,k)$-designs}

Let  $R$ be the reduced incidence matrix of a $3$-$(n,k)$-design. We now rewrite the properties from the previous subsection in terms of linear constraints. 

Since $k$ can be any integer, in the linear programming formulation considered below, $k$ is nonegative integer variable. Recall that the incidence matrix $A$ of a $3$-$(n,k)$-design is obtained from the reduced matrix $R$ by repeating some of its rows.

Let $z$ be a row vector with positive integer values indicating the number of times each row of matrix $R$ is repeated as a row in matrix $A$. We rewrite the property (P3) of incidence matrix $A$ (the sum in each column of matrix $A$ is $k$) as the following linear relation:

  \begin{equation}\label{LC1}(z|k)\left(%
\begin{array}{ccc}
&  &  \\
& R &  \\
&  &  \\
\hline
-1& .  &-1\\
\end{array}%
\right)=0.\end{equation}
Note that some values of the variables $z_i$ must be equal to $1$. Let $F$ (frozen) denote the set of indices  $i$  of variables $z$ (rows of $R$) such that
 the matrix $\left(%
\begin{array}{c}
R\\
R_i
\end{array}%
\right)$, obtained by appending of the $i$-th row of  $R$ do not fulfill $3$-orthogonality. Thus we obtain the following constraints:

\begin{equation}z_i=1, i\in F.\end{equation}

The set of indices $F$ can be obtained by a complete search of the rows of the matrix $A$ prior to formulating the linear programming problem. We write the constraints above in the form of an integer linear programming problem and  denote it as $\text{ILP}(R)$:

$$
\text{Integer variables}:$$ $$z_i, i \in \{1,\ldots, Nrows(R)\}, k,$$
$$\text{Subject to constraints:}$$
$$z_i\geq 1,  i \in \{1,\ldots, Nrows(R)\},k \geq 6,$$ 
$$z_i=1, i\in F,$$
$$(z|k)\left(%
\begin{array}{ccc}
&  &  \\
& R &  \\
&  &  \\
\hline
-1& .  &-1\\
\end{array}%
\right)=0.$$

Note that  the objective function is absent in this formulation. The constraint $k \geq 6$ is related to the fact that the previously solved cases $k = 4$ and $5$  are not considered. Specifically, the minimum number of blocks $d_{3,n,4}$ in a $3$-$(n,4)$-design   is $5$ and $14$, respectively, see Theorems \ref{T_1} and  \ref{T:t=3}.   

In the next paragraph, a description of the computer search program for reduced incidence matrices when $k$ is even is be provided, based on the problem described above. Note that the constraint on the parity of $k$ is not present in the formulation of $\mbox{ILP}(R)$, which indirectly shoes Theorem \ref{T:t=3} for odd $k$, given the results of the search (see the end of the next subsection).

\subsection{ILP-based search of reduced incidence matrices of  $3$-$(n,k)$-designs for even $k$}
By Theorems \ref{T:main}.2 and \ref{T:t2}.2 for even $k$, $n\geq k+3$, the minimum number of blocks  $d_{3,n,k}$ in a $3$-$(n,k)$-design is not less than  $7$:
$$d_{3,n,k}\geq d_{2,n-1,k-1}+1=7.$$

On the other hand, Construction \ref{Constr_2}, applied two times to a $1$-$(n-4,k-2)$-design with $3$ blocks (which exists for $n\geq k+4$), see (\ref{ubpund1}), gives the following inequality:

$$d_{3,n,k}\leq 12.$$

Therefore we have the following.
\begin{proposition}
   Let $k$ be even, $k\geq 6$ and $n\geq k+4$. Then we have that  
   $7\leq d_{3,n,k}\leq 12$.
\end{proposition}

In this paragraph, we  describe an algorithm for finding reduced incidence matrices of $3$-$(n,k)$-designs for even $k$, which is aimed to narrow the gap in the Proposition above.  Specifically, we would like  to show that $d_{3,n,k} = 12$ for all $n \geq k + 4$, except for a finite number of values for $k$. The algorithm is written in  MAGMA \cite{MAGMA} and is available in Appendix 1.

Firstly, we will describe the basic components of the program.

$1.$ Number of columns. The number of columns $d_{3,n,k}$ in the considered reduced matrices satisfies the inequality $7 \leq d_{3,n,k} \leq 11$. In the program, $d_{3,n,k}$ is an input parameter denoted as $nBlocks$, where $7 \leq nBlocks \leq 11$.

$2.$ Number of ones in rows. The number of ones in each row of the reduced incidence matrix of a $3$-$(n,k)$-design for even $k$, $k \geq 6$, is at least $6$, see $(P4')$. In the program, this constraint is represented via constant $minWeight = 6$.

$3.$ Iterative process. The program works iteratively   increasing the number of rows in the matrices within the set setOfMatrices. If a matrix (denoted by $mTemp$ in program) obtained from a matrix $M$ in the set setOfMatrices by appending a row satisfies certain necessary conditions (see more details in item 6 below) for the reduced incidence matrices of $3$-$(n,k)$-designs then the problem $\text{ILP}(mTemp)$ is launched. If $\text{ILP}(mTemp)$ has at least one solution, the matrix is output to the console.

$4.$ Isomorphism rejection. To reduce the search by excluding matrices obtained by permutations of rows and columns, pruning is done based on the canonical graph, which is implemented using the built-in function in MAGMA for computing the canonical graph $CanonicalGraph()$.

$5.$ Main program loop. At the start of the program, a  matrix consisting of one row is added to the set $setOfMatrices$. This row consists of $x$ consecutive ones, where $x \in \{minWeight, nBlocks\}$.

The program with the input parameter $nBlocks$ iterates over all $x$ in $\{minWeight, \ldots, nBlocks\}$. Each iteration is composed of several sub-iterations, each of which increases the number of rows in the matrices by $1$.

Every sub-iteration consists of the following steps:

5.1 Appending. For each matrix $M$ in the set $setOfMatrices$, row-vectors (different from the rows of matrix $M$) are appended such that $3$-orthogonality holds of the matrix M with appended vector. Without restriction of generality, we take the appended vector having the weight not less than that of the last row of $M$. Note that when adding a row to a matrix consisting of only one row, the check for $3$-orthogonality is not performed because at this stage there is less than three rows in matrix. 

5.2 Usage of isomorphism rejection. For each matrix $mTemp$ with an appended row, its full invariant is computed based on the canonical graph. If the canonical graph is new (i.e., it is not in the set of graphs stored in $setOfCanonicGraphs$), it is added to the list of canonical graphs and the matrix is added to the set $setOfMatricesForNextIteration$. At the end of the sub-iteration, when all possible vectors have been tested for all matrices for appending, the variable set $setOfMatrices$ is assigned to $setOfMatricesForNextIteration$ and the next sub-iteration starts with the updated set of matrices. A sub-iteration ends when no vectors can be added to the matrices in the set $setOfMatrices$ without violating $3$-orthogonality.

5.3. Conditions for launching  problem $ILP(mTemp)$. Inside the main loop, for each matrix $mTemp$ where the columns are pairwise distinct, $3$-orthogonality holds and the canonical graph is new for the set $setOfCanonicGraphs$ the linear programming problem $\mbox{ILP}(mTemp)$ is solved.

{\bf The results the search.}

The algorithm for fixed values of $nBlocks$ from $7$ to $10$ is relatively fast and fits within the two-minute limit for computations on the online Magma calculator \cite{MAGMAcalc}. For $nBlocks=11$ it runs for about 10 minutes (which requires a licensed version of MAGMA).

By running program for all $nBlocks$ in range from $7$ to $11$, the following matrices were obtained: $J_7-E_7$, $J_9-E_9$, $J_{11}-E_{11}$, where $J_l$ and $E_l$ are all-ones and identity matrices of size $l$ respectively. It is easy to see that these matrices are the incidence matrices of the $3$-$(7,6)$, $3$-$(9,8)$, and the $3$-$(11,10)$-designs obtained using Construction \ref{Constr_1}. For  any of the matrices $J_7-E_7$, $J_9-E_9$, $J_{11}-E_{11}$ consider two different rows and take a  copy of one of them. We see that three such  rows have exactly $5$, $7$ and $9$ respectively common ones.
So any of the matrices 
 $J_7-E_7$, $J_9-E_9$, $J_{11}-E_{11}$ can not be extended by adding duplicates of any of its rows and fulfill $3$-orthogonality. 
In other words, there are no other incidence matrices of $3$-designs that have these three matrices as their the reduced incidence matrices, which are in turn incident matrices as well. We conclude the following.

\begin{theorem}\label{S:t3even}
1.\cite{Pot}\cite{WZ} For any $n\geq 5$ we have 
$d_{3,n,4}=5$.

2. If $k$ is even, $k\geq 6$ and $n\geq k+4$ then $$d_{3,n,k}=\begin{cases}
     7, k=6,\\
     9, k=8,\\
     11, k=10,\\
     12, k\geq 12
 \end{cases}.$$

\end{theorem}

We summarize the results of the minimum distance study of codes $C_{t,n,k}$ in Table \ref{Tr1}.

\begin{table} \caption{\normalsize Minimum distances $d_{t,n,k}$ for   large enough $n$}\label{Tr1}\small

\begin{center}
\noindent
\hskip-1.31cm
\begin{tabular}{|c|c|c|c|c|c|c|c|c|c|c|c|c|}
  \hline
  $t/k$ &$2$&$3$&$4$& $5$& $6$  & $7$ & $8$ & 9  &10 &11 &12 & 13\\
   \hline
 $4$ & $-$ & $-$ &  $-$& $6^1$         & $21^8\ldots  28^9$ & $8^{7,10}$ & $32^3$ & $10^{7,10}$ & $25^{8}$\ldots $28^9$ &$12^{7,10}$ &  $32^3$ & $13^7 \ldots 14^{10}$\\
 $3$ & $-$ & $-$ &  $5^1$& $14^6$ &   $7^5$          & $16^6$ &$9^5$ &$14^6$ & $11^5$&  $16^6$ & $12^5$ & $14^6$  \\
 $2$ & $-$ & $4^1$ & $7^4$ & $6^4$ & $8^4$ & $6^4$ & $7^4$ & $6^4$ & $8^4$ & $6^4$  & $7^4$ & $6^4$ \\
 $1$ & $3^1$ & $4^2$ & $3^2$ & $4^2$ & $3^2$ &$4^2$& $3^2$& $4^2$ & $3^2$ & $4^2$ & $3^2$ & $4^2$   \\
  \hline
 \end{tabular}

 \end{center}

{\it $^1:$ Theorem \ref{T_1} \cite{Pot} \cite{WZ}: $d_{k-1,n,k}=k+1$} 

{\it {$^2$: \mbox{ Proposition } \ref{p2} } }

{\it {$^3$: \mbox{ Theorem } \ref{T2}} } 

{\it {$^4$: \mbox{ Theorem } \ref{T:t2}} }

{\it {$^5$: \mbox{ Theorem } \ref{T:t=3}} }

{\it {$^6$: \mbox{ Theorem }\ref{S:t3even} } }

{\it {$^7$: \mbox{Lower bound from Theorem \ref{T:main}.2} }} 

{\it {$^8$: \mbox{Lower bound from Theorem \ref{T:main}.4} }}

{\it {$^9$: \mbox{ Upper bound }(\ref{ubpund1}) } }

{\it {$^{10}$: \mbox{ Upper bound (\ref{Ub_constr_1}) for small $k$   } }}

 \end{table}

\section{Quasicyclic LDPC codes from Wilson inclusion matrices}\label{S_QC}

We developed a \Csharp\, program that  constructs the quasicyclic codes avoiding cycles of length four or six  (i.e. exponent matrices fulfilling  (\ref{C4}) or (\ref{C6})) given the base matrix. The idea behind the construction was as follows. Firstly, we purged all cycles incident to the some (randomly chosen but not all) variable vertices in a random order 
 incident to $4$- and $6$- cycles. We used randomness to avoid oscillations in this stage.  Then we purged the   cycles incident to the remaining variable vertices, arranged in decreasing number of incident $4$- and $6$- cycles. 

 We compared the obtained quasicyclic LDPC codes from Wilson-type exponent matrices against random codes from \cite{MacKay} having the same sizes parity check matrices and column sums. We produced the following three codes from our \Csharp\, program:

$\bullet$ A code of length $2310$ and dimension $1816$ obtained from the matrix $W_{2,10,4}$ by circulant lifting with $qc=11$ and no cycles of length four in Tanner graph.

$\bullet$ A code of length $27720$ and dimension $19801$ obtained from the matrix $W_{2,9,4}$ by circulant lifting with $qc=220$ and no cycles of length four in Tanner graph.

$\bullet$ A code of length $27720$ and dimension $19801$ obtained from the matrix $W_{2,9,4}$ by circulant lifting with $qc=220$, without cycles of lengths four and six in Tanner graph.
All obtained codes have a prefull rank, i.e. one less than the number of their rows.

The obtained codes were tested for the following decoders and channels: min-sum algorithm with layered schedule \cite{minsum} and AWGN channel;  multiple gradient descent bit flipping  \cite{GDBF} and binary symmetric channel, see Fig. $1$ and Fig. $2$. From the plots we see that both codes cycle four free and cycle four and six free codes from Wilson matrices behave similarly to the codes  from MacKay matrices.

\begin{figure}[H] 
\centering
\includegraphics[width=1\textwidth]{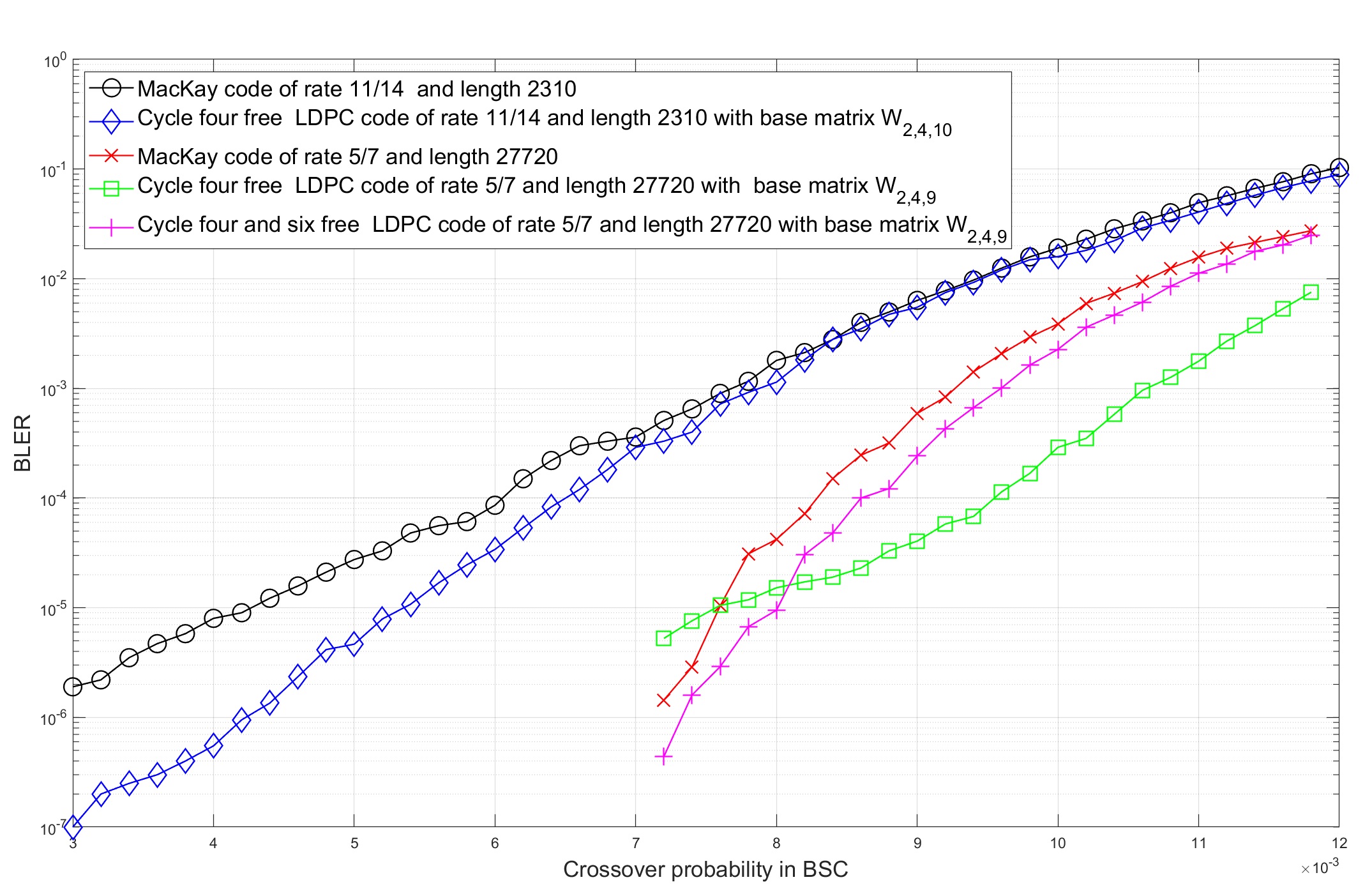} 
\caption{Decoding performance of 
    QC LDPC codes from Wilson matrices
	vs MakKay type 1A code 
    under Multi GBDF decoder
  with 30 iterations max in binary symmetric channel.}
\label{fig:first_graph}
\end{figure}

 \begin{figure}[H] 
\centering
\includegraphics[width=1\textwidth]{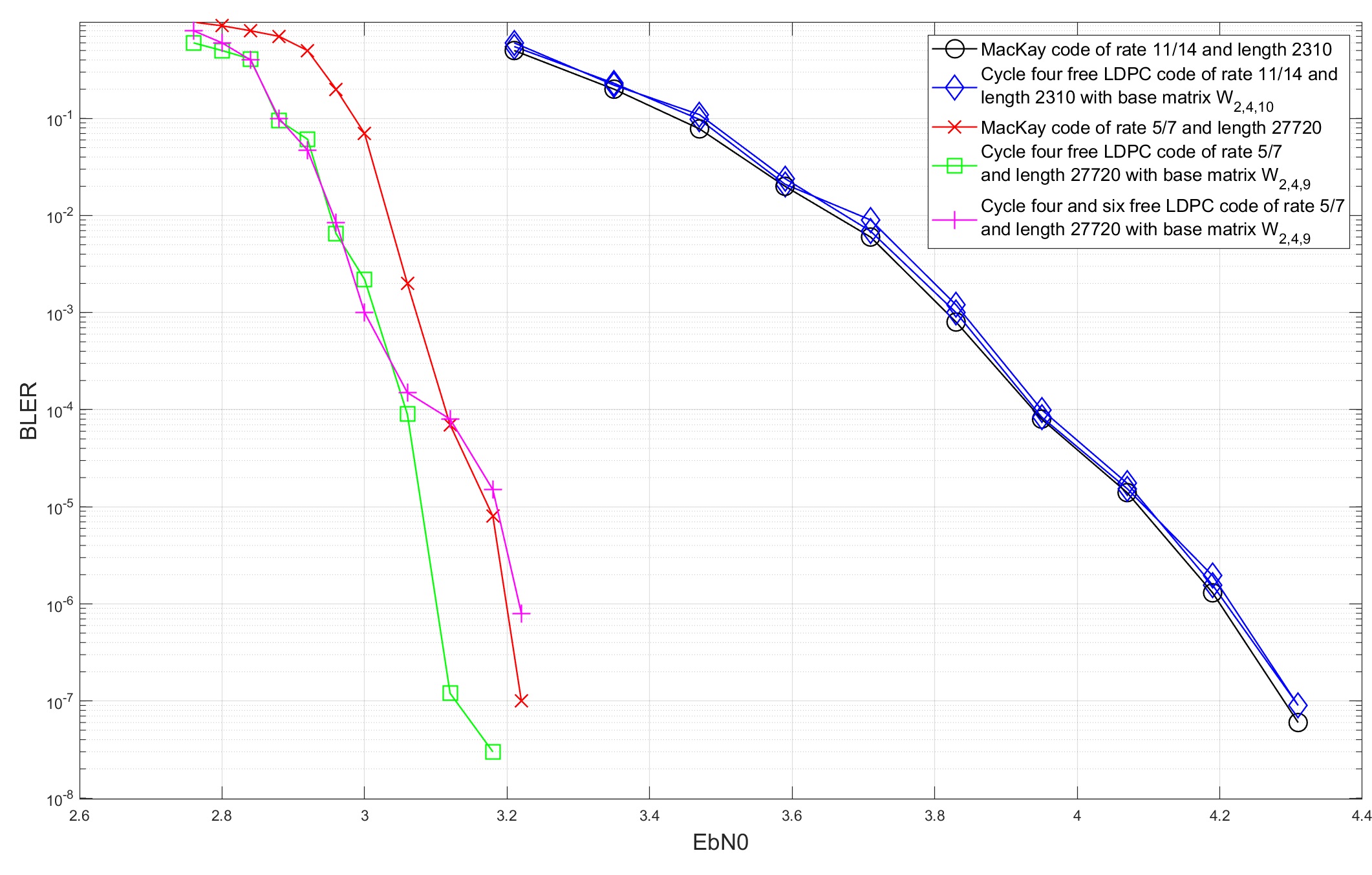} 
\caption{Decoding performance of 
    QC LDPC codes from Wilson matrices
	vs MakKay type 1A code 
    under layered min sum algorithm 
  with 30 iterations max in AWGN channel.}
\label{fig:first_graph}
\end{figure}

{\bf Acknowledgements.} The authors would like to express their gratitude to Alexey Frolov for a talk on locally recoverable codes at the online seminar 'Coding theory', Vladimir Potapov for talks on unitrades at the seminar '2024-ary quasigroups and related topics'. These contributions significantly directed their focus towards the current research.  The authors are profoundly thankful to Evgeny Vdovin for  providing the licensed MAGMA software, which was indispensable for this work.

\appendix
\section{MAGMA program  for finding binary $3$-$(n,k)$-designs with small number of blocks}
\begin{verbatim}
nBlocks:=9;
minWeight:=6;

V0:=VectorSpace(GF(2),nBlocks);

isOrthogonal3:= function(M,v);
	numeration:= Subsets({1..Nrows(M)},2);
	for num in numeration do
		numList:=[];
		for id in num do
			Append(~numList, id);
		end for;
		countOf3:=0;
		for j in [1..Ncols(M)] do
			if (M[numList[1]][j] eq 1) and (M[numList[2]][j] eq 1) and (v[j] eq 1) then
				countOf3:= countOf3 + 1;
			end if;
		end for;
		if (countOf3) mod 2 ne 0 then
			return false;
		end if;
	end for;
	return true;
end function;

isRowCanBeReused:= function(M,rowID);
	for row in [1..Nrows(M)-1] do
		countOf3:=0;
		for j in [1..Ncols(M)] do
			if (M[row][j] eq 1) and (M[rowID][j] eq 1) and (M[rowID][j] eq 1) then
				countOf3:= countOf3 + 1;
			end if;
		end for;
		if (countOf3) mod 2 ne 0 then
			return false;
		end if;
	end for;
	return true;
end function;

findLinearIntegerSolutions:=procedure(mTemp);
    mTempJoint:= 
VerticalJoin(Matrix(Integers(),mTemp),Matrix(Integers(),1,Ncols(mTemp),[(-1)^^Nco
ls(mTemp)]));
	nRowsOfM:=Nrows(mTempJoint);
	nColsOfM:=Ncols(mTempJoint);
	L:=LPProcess(Integers(), nRowsOfM);
	R:=Matrix(Integers(),nColsOfM,1,[0^^nColsOfM]);
	MjT:=Transpose(mTempJoint);
	possibleNumberOfUse:=Matrix(Integers(),nRowsOfM,nRowsOfM,[0^^nRowsOfM*nRowsOfM]);

	for i in [1..nRowsOfM-1] do
		if isRowCanBeReused(mTempJoint, i) then
			possibleNumberOfUse[i][i]:=0; 
		else
			possibleNumberOfUse[i][i]:=1; 
		end if;
	end for;					
					
	maxPossibleUseOfRows:=[1^^nRowsOfM];


	AddConstraints(L,MjT,R: Rel:="eq");
	
AddConstraints(L,ScalarMatrix(nRowsOfM,Integers()!1),Matrix(Integers(),nRowsOfM,1,[1^^(nRowsOfM)]) : Rel:="ge");
        AddConstraints(L,Matrix(Integers(), 1,nRowsOfM,[0^^(nRowsOfM-1)]cat [ 1]),Matrix(Integers(), 1,1,[ 6]): Rel:="ge");

	X:=Solution(L);
	if {X[1][i] eq 0:i in [1..Ncols(X)]} ne {true} then
		"Matrix:"; mTempJoint;
		"Solution:"; X;
		"---------------------";
	end if;
end procedure;

convertMatrixToCanonicGraph:= function(M)
	Z0:=ZeroMatrix(Integers(),Ncols(M),Ncols(M));
	Z1:=ZeroMatrix(Integers(),Nrows(M),Nrows(M));
	M1:=VerticalJoin(HorizontalJoin(Z0,Transpose(M)),HorizontalJoin(M,Z1));
	Can:=CanonicalGraph(Graph<Nrows(M1)|M1>);
	return Can;
end function;

isColumsDifferent:=function(M);
	MT:=Transpose(M);
	for row1 in [1..Nrows(MT)-1] do
		for row2 in [row1+1..Nrows(MT)] do
			if MT[row1] eq MT[row2] then
				return false;
			end if;
		end for;
	end for;
	return true;
end function;

time for countOfOnes in [minWeight..nBlocks] do
	x:=[1^^countOfOnes]cat [0^^(nBlocks-countOfOnes)];
	printf "Starting for countOfOnes = %o\n", countOfOnes; 
	setOfMatrices:={Matrix(GF(2),1,nBlocks,x)};
		
	V:={v:v in V0|Weight(v) ge countOfOnes};
	repeat
	
		setOfCanonicGraphs:={};
		setOfMatricesForNextIteration:={};
		repeat 
			M:=Random(setOfMatrices);
			Exclude(~setOfMatrices,M);
			for v in V do
				if Weight(v) ge Weight(M[Nrows(M)]) then  
					checkVinM:= 0;
					for j in [1..Nrows(M)] do
						if v eq M[j] then
							checkVinM:=1;
							break j;	
						end if;
					end for;
	
					if((checkVinM eq 0) and ((Nrows(M) le 1) or isOrthogonal3(M,v)))  then 
						mTemp:= VerticalJoin(M,Matrix(GF(2),1,nBlocks,Eltseq(v)));
						mTempGraph:=convertMatrixToCanonicGraph(Matrix(Integers(),mTemp));
						if mTempGraph  notin setOfCanonicGraphs then
       	                                       		if isColumsDifferent(mTemp)  then
													findLinearIntegerSolutions(Matrix(Integers(), mTemp));
							end if;
							Include(~setOfCanonicGraphs,mTempGraph);
							Include(~setOfMatricesForNextIteration,mTemp);
						end if;
					end if;
				end if;
			end for;
		until #setOfMatrices eq 0;
		setOfMatrices:= setOfMatricesForNextIteration;
	until #setOfMatricesForNextIteration eq 0;
	printf "Finished for countOfOnes = %o\n", countOfOnes; 
end for;
"Program finished.";
\end{verbatim}

\end{document}